\documentclass[12pt]{article}
\usepackage{amsmath,amsfonts,amsthm,amscd,amssymb,graphicx}

\usepackage{graphicx}

\usepackage{amssymb}
\usepackage[colorlinks,citecolor=blue,linkcolor=blue]{hyperref}
\usepackage{xcolor}

\def\cit{{\Bbb C}}

\def\P{{\Bbb P}}
\def\eps{\varepsilon}

\newtheorem{e-proposition}[theorem]{Proposition}

\newtheorem{e-definition}[theorem]{Definition\rm}

\newtheorem{result}{Result}[section]

\def\og{\leavevmode\raise.3ex\hbox{$\scriptscriptstyle\langle\!\langle$~}}
\def\fg{\leavevmode\raise.3ex\hbox{~$\!\scriptscriptstyle\,\rangle\!\rangle$}}

\def\beq{\begin{equation}}
\def\eeq{\end{equation}}

\begin{document}



%

\centerline{\Large \bf Onset of nonlinear instabilities}

\bigskip

\centerline{\Large \bf in monotonic viscous boundary layers}

\bigskip

\centerline{D. Bian\footnote{Beijing Institute of Technology, School of Mathematics and Statistics, Beijing, China.
Email: biandongfen@bit.edu.cn}, 
E. Grenier\footnote{UMPA, CNRS  UMR $5669$, Ecole Normale Sup\'erieure de Lyon, Lyon, France.
Email: Emmanuel.Grenier@ens-lyon.fr}}

. 


\subsubsection*{Abstract}


In this paper we study the nonlinear stability of a shear layer profile for Navier Stokes equations near a boundary.
More precisely, we investigate the effect of cubic interactions on the growth of the linear instability.
In the case of the exponential profile  we obtain that the nonlinearity tames the linear instability.
We thus conjecture that small perturbations grow until they reach a magnitude $O(\nu^{1/4})$ only, forming small rolls
in the critical layer near the boundary.
The mathematical proof of this conjecture is open.


\section{Introduction}


In this paper we consider the incompressible Navier Stokes equations in an half plane
 \beq \label{NS1} 
\partial_t u^\nu + (u^\nu \cdot \nabla) u^\nu - \nu \Delta u^\nu + \nabla p^\nu = f^\nu,
\eeq
\beq \label{NS2}
\nabla \cdot u^\nu = 0 ,
\eeq
together with the Dirichlet boundary condition
\beq \label{NS3} 
u^\nu = 0 \qquad \hbox{for} \qquad y = 0.
\eeq
We are interested in the stability of a shear layer profile
$$
U(y) = (U_s(y),0).
$$ 
Note that this shear layer profile is a stationary solution of Navier Stokes equations
provided we add the forcing term $f^\nu = (- \nu \Delta U_s,0)$. We assume that $U_s$ is a smooth function, with
$U_s(0) = 0$, $\partial_y U_s(0) \ne 0$ and that $U_s$ converges at $+\infty$ to some constant $U_+ \ne 0$.
Note in particular that this includes profiles like the exponential profile $U_s(y) = 1 - \exp( - \delta y)$
where $\delta$ is a positive constant,  but not Couette's profile.

The stability of such shear layer profiles is one of the most classical questions in fluid dynamics and
has been intensively studied in physics since the beginning of the twentieth century, in particular by Rayleigh,
Prandtl, Orr, Sommerfeld, Tollimien, Schlichting, C.C. Lin. 
More recently, the rigorous study of the linear instability of shear layers has been pioneered by \cite{GGN3}.

This question is also closely related to the study of the inviscid limit of Navier Stokes equations near a boundary.
For initial data with analytic regularity, it is known that, in small times, solutions of Navier Stokes equations converge to solutions 
of Euler equations, up to adding Prandtl's boundary layers \cite{SammartinoCaflisch1}, \cite{SammartinoCaflisch2}. 
This work has been then extended  for instance in \cite{GVM}, \cite{Mae}. For initial data with Sobolev regularity,
this question is widely open. In \cite{GN2}, the authors construct a sequence of smooth solutions of
the Navier Stokes equations which can not be described as solutions of Euler with Prandtl boundary layer in 
$L^\infty$. This example is precisely a boundary layer profile, whose instability in $L^\infty$ is completely described.
A natural question is thus to investigate whether {\it all} boundary layer profiles are unstable in $L^\infty$.
The aim of this article is to initiate a formal analysis of this question.

It turns out that non trivial shear layer profiles are always linearly unstable with respect to Navier Stokes equations \cite{GGN3}.
Two classes of linear instabilities appear.

\begin{itemize}

\item "Inviscid instabilities": instabilities which persist as $\nu$ goes to $0$. According to Rayleigh's criterium,
such instabilities only occur for profiles
$U_s$ with inflection points. They exhibit scales in $t$, $x$ and $y$ of order $O(1)$.
The corresponding eigenvalue $\lambda$ has a real part $\Re \lambda$ of order $O(1)$.

\item "Long wave instabilities" \cite{BG1}: these instabilities arise even in the case of concave profiles $U_s$, such that
$U_s'' < 0$. They do not persist as $\nu$ goes to $0$ and 
are characterized by a strong spatial anisotropy since their sizes are of order $O(1)$ in $y$ but of order
$O(\nu^{-1/4})$ in $x$. Moreover they grow very slowly, within time scales of order $O(\nu^{-1/2})$.
The corresponding eigenvalue $\lambda$ has a real part $\Re \lambda$ of order $O(\nu^{1/2})$.

\end{itemize}

It has been proven in \cite{GN2} that linear "inviscid instabilities" lead to a nonlinear instability in $L^\infty$ norm.
More precisely, there exists arbitrarily small perturbations of $U$ (small like $\nu^N$ in $H^N$ for arbitrarily large
$N$), which lead to a $O(1)$ difference on the nonlinear solution in $L^\infty$ within time scales of order $O(\log \nu^{-1})$.
This gives an example of sequence of solutions for which Prandtl's analysis is not true and proves that
the classical works of R.E. Caflisch and M. Sammartino
\cite{SammartinoCaflisch1}, \cite{SammartinoCaflisch2} can not be extended to solutions with Sobolev regularity.

Such a result remains unknown in the case of "long wave instabilities". 
We know that concave profiles $U_s$ are linearly unstable, but we do not know, up to now, whether they are nonlinearly unstable.
It is not known whether shear layers develop nonlinear instabilities which reach order $O(1)$ in $L^\infty$.

The main difficulty is that the instability grows too slowly, and hence quadratic terms can not be neglected.
Let us illustrate this statement on a simple ordinary differential equation and consider the classical bifurcation model equation
\beq \label{model}
\dot \phi = \lambda \phi + A |\phi|^2 \phi + O(\phi^5)
\eeq
where $\phi$ is a scalar, $\Re \lambda > 0$ is a small parameter and $A \in \cit$.
The evolution of a small perturbation of $\phi = 0$ depends on the sign of $\Re A$. 
If $\Re A < 0$ then any small perturbation grows until it reaches
a size $O(\lambda^{1/2})$.
If $\Re A > 0$ then any small perturbation grows and reaches $O(1)$. 

In this paper we are in a similar situation. The most unstable linear mode  grows like $e^{\lambda t}$ with
$\Re \lambda = O(\nu^{1/2})$. More precisely it is of the form
$$
v_{lin} = \eps e^{\lambda t} \nabla^\perp \Bigl( e^{i \alpha x}  \psi_{lin}(y) \Bigr) + c.c.
$$
where $e^{i \alpha x} \psi_{lin}(y)$ is the corresponding stream function, $\eps$ is a small parameter (initial size of the
perturbation), $\alpha$ is the horizontal wave number, of order $O(\nu^{1/4})$, and where
$c.c.$ denotes the complex conjugate.

Starting from this linear perturbation it is classical to construct an approximate solution of Navier Stokes equations in the form of
$$
u^{app} = U + v_{lin} + v_{quad} + v_{cubic} + ...
$$
where $v_{quad}$ gathers quadratic terms and $v_{cubic}$ cubic ones.
The proof of \cite{GN2} relies on the construction of such an approximate solution at any order
$$
u^{app} \sim U +  \sum_{n \ge 0} \eps^n e^{n \lambda t} u_n 
$$
and on the proof of the convergence of the corresponding infinite series. 
Such an approach however is impossible in our case since, as will be clear in the sequel, $u_{lin}$ is not larger than the other terms.

The aim of this paper is to study these quadratic and cubic terms. We first observe that, as $v_{lin}$ has horizontal wavenumbers
$\alpha$ and $-\alpha$, $v_{quad}$ has horizontal wavenumbers $2\alpha$ and $- 2 \alpha$, and $v_{cubic}$ has
wavenumbers $\alpha$, $-\alpha$, $3 \alpha$ and $- 3 \alpha$. Thus the first feedback of the nonlinearity on the linear instability
occurs with $v_{cubic}$. A natural question is to know whether this cubic term enhances or tames the linear one.
More precisely: what is the sign of the projection of the horizontal wavenumber $+ \alpha$ of $v_{cubic}$ on the
horizontal wavenumber  $+ \alpha$ of $v_{lin}$ ? 

Let us denote this projection by
$$
\P(v_{cubic},v_{lin},+\alpha) .
$$
If $\Re \P(v_{cubic},v_{lin},+\alpha) > 0$ then we expect that the nonlinear interactions enhance the linear instability.
In this case, small perturbations would grow and reach a size of order $O(1)$ in $L^\infty$.
On the contrary, if $\Re \P(v_{cubic},v_{lin},+\alpha) < 0$, then we expect nonlinear interactions
to tame the instability, which is then likely to saturate when $\lambda v_{lin}$ and $v_{cubic}$ are of the same order, namely
of order $O(\nu^{1/4})$. The situation is very close to the "bifurcation" scenario discussed on the model equation (\ref{model}).

The main result of the paper is the following

\begin{result}
For exponential boundary layer profiles of the form 
$$
U_s(y) = 1 - e^{-\delta y},
$$
numerical computations show that
$$
\Re \P(v_{cubic},v_{lin},+\alpha) < 0.
$$
\end{result}

This result relies on careful numerical and formal computations which are detailed in the forthcoming sections.
The computations are delicate since the solution has three spatial scales, namely
$O(1)$ (that of the shear layer itself), $O(\nu^{1/4})$ (size of the so called "critical layer", see \cite{GGN3}), and
$O(\nu^{-1/4})$ (horizontal instability size, and also recirculation size). 
Moreover the vorticity is concentrated in the critical layer of size $O(\nu^{1/4})$ and has a singular behavior as $\nu \to 0$. 
To handle these difficulties we follow a "mixed" approach. We  describe the solution
by series in $y$ and $y \log y$ near the boundary, by a numerical grid in the critical layer and by another numerical grid away from the
boundary, which allow to deal with arbitrarily small viscosity.
  
\medskip

Let us now discuss the consequences of such a result for the nonlinear instability of shear layers.
The end of this section is purely conjectural. The situation is very close to that of a bifurcation.
The role of the bifurcation parameter is played by $\lambda = O(\nu^{1/2})$. The bifurcation occurs at $\nu = 0$.
Generically when $\nu > 0$, two eigenvalues $\lambda$ and $\bar \lambda$ emerge from $\Re \lambda = 0$, 
exactly as in an Hopf bifurcation.

This leads to linearly growing modes, with very slow growths, of size $\Re \lambda \sim O(\nu^{1/2})$.
Quadratic, cubic and higher order interactions take place, which can no longer be neglected when they
reach the size of $\lambda u_{lin}$, namely when the linear instability reaches a magnitude $O(\nu^{1/4})$.
The instability may then saturate at the magnitude $O(\nu^{1/4})$, as is the case in the classical Hopf bifurcation scenario
\cite{Haragus}.

We recall that the domain of instability in $\alpha$ is $C_1 \nu^{1/4} \le \alpha \le C_2 \nu^{1/6}$ for some constants
$C_0$ and $C_1$ \cite{GGN3}. If $2 \alpha \le C_2 \nu^{1/6}$ then $u_{quad}$, which has horizontal wavenumbers
$\pm 2 \alpha$, also creates linear instabilities which grow and saturate, leading to a possible cascade of instabilities.
If $2 \alpha > C_2 \nu^{1/6}$, then wavenumbers $\pm 2 \alpha$ are linearly stable and the situation is close to Rayleigh Taylor
instability. In this case we conjecture the existence of solutions which converge to "rolls" as $t$ goes to infinity.
However, in strong contrast with the usual assumptions in bifurcation theory, there is no "spectral gap" in our case.

The paper is constructed as follows. In part \ref{part2} we detail the construction of an approximate solution up to the third order.
Part \ref{part3} is devoted to the design of a general strategy to invert Orr Sommerfeld equations. 
The growing linear modes  and the corresponding
modes of the adjoint of Orr Sommerfeld equations are constructed in part \ref{linear}.
Part \ref{part5} is devoted to some details on the computations in the particular case of an exponential profile.


\section{Principle of the construction \label{part2}}


 
 \subsection{Orr Sommerfeld equations}
 
 
 Let us first introduce the classical Orr Sommerfeld equations.
  We refer to \cite{GN} for more details on all these aspects.
Let $L$ be the linearized Navier Stokes operator near the shear layer profile $U$, namely
\beq \label{linearNS}
L v  = (U \cdot \nabla) v + (v\cdot \nabla) U - \nu \Delta v + \nabla q,
\eeq
with $\nabla \cdot v = 0$ and Dirichlet boundary condition. We want to study the resolvant of $L$, namely to
study the equation
\beq \label{resolvant}
(L + \lambda) v = f,
\eeq
where $f$ is a given forcing term and $\lambda$ a complex number.
Taking advantage of the divergence free condition, we introduce the stream function $\psi$ and take the Fourier
transform in $x$ and the Laplace transform in $t$, which leads to look for solutions of the form
$$
v = \nabla^\perp \Bigl( e^{i \alpha (x - c t) } \psi(y) \Bigr) .
$$
Note that $\lambda = - i \alpha c$. We also take the Fourier and Laplace transform of the forcing term $f$
$$
f = \Bigl( f_1(y),f_2(y) \Bigr) e^{i \alpha (x - c t) } .
$$
Taking the curl of (\ref{resolvant}) we then get the classical Orr Sommerfeld equations
\beq \label{Orrmod0}
Orr_{\lambda,\alpha,\nu}(\psi) =  (U_s - c)  (\partial_y^2 - \alpha^2) \psi - U_s''  \psi  
- { \nu \over i \alpha}   (\partial_y^2 - \alpha^2)^2 \psi =  i {\nabla \times f \over \alpha}
\eeq
where
$$
\nabla \times (f_1,f_2) = i \alpha f_2 - \partial_z f_1.
$$
The Dirichlet boundary condition
gives 
\beq \label{condOrr}
\psi(0) = \partial_z \psi(0) = 0.
\eeq
Let us define
\beq \label{epsilon}
\eps = {\nu \over i \alpha} .
\eeq
As $\nu$ goes to $0$, the Orr Sommerfeld operator degenerates into the Rayleigh operator
\beq \label{Rayleigh}
Ray_{c,\alpha}(\psi) = (U_s - c + \eps \alpha^2)  (\partial_y^2 - \alpha^2) \psi 
 - U_s''  \psi,
 \eeq
 which is a second order operator, together with the boundary condition 
 \beq \label{condRay}
 \psi(0) = 0.
 \eeq
 We have
 $$
 Orr_{\lambda,\alpha,\nu} = Ray_{c,\alpha} + Diff
 $$
 where
 $$
 Diff = - \eps (\partial_y^2 - \alpha^2) \partial_y^2 .
 $$
 Rayleigh operator is a good approximation of Orr Sommerfeld operator away from the boundary $y = 0$.
 Note that (\ref{Rayleigh}) degenerates at the critical layer $y_c$ defined by
 $$
 U_s(y_c) = c - \eps \alpha^2 .
 $$
In this critical layer we approximate the Orr Sommerfeld operator by the so called modified Airy operator ${\cal A}$
defined by 
$$
{\cal A} = Airy \, (\partial_y^2 - \alpha^2) 
$$
where
$$
Airy = (U_s - c + \eps \alpha^2) - \eps \partial_y^2.
$$
We have
$$
Orr_{\lambda,\alpha,\nu} = {\cal A} + ErrAiry
$$
where
$$
ErrAiry = - U''_s .
$$
Note that for very large $y$, solutions of $Ray_{c,\alpha}(\psi) = 0$ behave exponentially, like $e^{\pm \alpha y}$.

The analysis of Orr Sommerfeld equation has been detailed in \cite{GN} and we just recall a few basic facts here.
First, note that $\omega = -\Delta_\alpha \psi$ corresponds to the vorticity of the stream function $\psi$, where
$\Delta_\alpha = \partial_y^2 - \alpha^2$ and that (\ref{Orrmod0}) may be rewritten
\beq \label{Orrmod00}
	-(U_s - c) \omega - U''_s \psi + \eps \Delta_\alpha \omega 
	= i \alpha^{-1} \nabla \times f .
	\eeq
As $y$ goes to $+ \infty$, we note that $U''_s$ decays exponentially fast, like $e^{- \delta y}$, thus
$\omega$ decays at the same speed. On the contrary, $\psi$ may decay as slowly as $e^{- \alpha y}$.

 
 \subsection{Construction of an approximate solution}
 

The idea is to construct an approximate solution up to the third order, namely of the form
$$
u_{app} = U + \nu^N e^{\lambda t} u_{lin} + \nu^{2N} e^{2 \lambda t} u_q + \nu^{3N} e^{3 \lambda t} u_c +
c.c.
$$
($c.c.$ meaning complex conjugate), and to understand the geometry of $u_c$ with respect to $u_{lin}$, in order
to know whether cubic terms enhance or stop the linear instability.
Let us detail the computations to be done.

The starting point is a linear instability of Navier Stokes equations, namely an exponentially growing mode
of the corresponding linearized equations. 
Following the classical analysis we take advantage of the incompressibility condition to introduce the stream function,
take its Fourier transform in $x$ and its Laplace transform in time, and thus look for instable modes of the form
$$
u_{lin}^+ = \nabla^\perp \Bigl( e^{i \alpha (x - c t) } \psi_{lin}^+(y) \Bigr).
$$
Such vector fields are solutions of linearized Navier Stokes equations provided $\psi_{lin}$ satisfies the Orr
Sommerfeld equation
$$
Orr_{\lambda,\alpha,\nu}(\psi_{lin}^+) = 0.
$$
The spectral analysis of Orr Sommerfeld equation is recalled in section \ref{linear}
where it is proved that there exists unstable eigenmodes which are of the form
\beq \label{psilin}
\psi_{lin}^+ = \psi_{s,-} + a \psi_{f,-} \sim (U_s - c) e^{-\alpha y} + a Ai \Bigl( 2,\gamma (y - y_c) \Bigr),
\eeq
where $a = O(\nu^{1/4})$, $\gamma = O(\nu^{-1/4})$ and $\Im c = O(\nu^{1/4})$.
Splitting into horizontal wavenumbers $\alpha$ and $-\alpha$, we write 
$$
\psi_{lin} = \psi_{lin}^+ + \psi_{lin}^-
$$
and similarly for $u_{lin}$ and $\omega_{lin}$, $\psi_{lin}^-$ being the complex conjugate of $\psi_{lin}^+$. 
Note that the exponential in factor of $\psi_{lin}^-$ is $\exp(- i \alpha x + \bar \lambda t)$.
The associated velocity field and vorticity are
$$
u_{lin}^+ = \nabla^\perp \psi_{lin}^+ 
\sim \Bigl( \begin{array}{c}
 U_s' e^{- \alpha y} +a \gamma Ai(1,\gamma (y - y_c)) \cr
-i \alpha U_s e^{- \alpha y} - i \alpha a Ai(2,\gamma (y - y_c)) \cr
\end{array} \Bigr)
$$
and  
$$
\omega_{lin}^+ = \nabla \times u_{lin}^+ = -(\partial_y^2 - \alpha^2) \psi_{lin}^+
\sim -U_s'' e^{- \alpha y} - \gamma^2 a Ai(\gamma (y - y_c)) .
$$
  The next term of the expansion is $u_q$, solution of
 \beq \label{uq1}
 \partial_t \omega_q + (U \cdot \nabla) \omega_q - \nu \Delta \omega_q + U_s'' \partial_x \psi_q  = Q_1,
 \eeq
where $\omega_q = \nabla \times u_q$
and 
$$
Q_1 = -  (u_{lin} \cdot \nabla) \omega_{lin} 
=  -  (u_{lin}^+ \cdot \nabla) \omega_{lin}^+   -  (u_{lin}^- \cdot \nabla) \omega_{lin}^-  ,
$$ 
since $+ / -$ and $- / +$ interactions vanish.
Note that $Q_1$ is the sum of two terms, one with wavenumber $2 \alpha$ and the other with wavenumber $- 2 \alpha$,
hence
$$
Q_1 = Q_1^+ + Q_1^-.
$$
Using Orr Sommerfeld equation we get that
$$
\psi_q = \psi_q^+ + \psi_q^-
$$
where
$$
Orr_{2 \lambda, 2 \alpha,\nu}(\psi_q^+) =   {i \over 2 \alpha} Q_1^+ .
$$
Note that we also need to enforce the boundary conditions 
$$
\psi_q^\pm(0) = \partial_y \psi_q^\pm(0) = 0.
$$
Note also that the exponentials in factor of $\psi_q^\pm$ are $\exp( 2 i \alpha x + 2 \lambda t)$
and $\exp( - 2 i \alpha x + 2 \bar \lambda t)$. 
The third term of the expansion is $u_c$, solution of
\beq \label{uc1}
\partial_t \omega_c + (U \cdot \nabla) \omega_c - \nu \Delta \omega_c + U_s'' \partial_x \psi_c  = Q_2,
\eeq 
where
$$
Q_2 = - (u_{lin} \cdot\nabla) \omega_q - (u_q \cdot\nabla) \omega_{lin}. 
$$
Note that $Q_2$ is a sum of four terms, with wavenumbers $- 3 \alpha$, $- \alpha$, $\alpha$, $3 \alpha$, namely
$$
Q_2 = Q_2^{-3} + Q_2^{-1} + Q_2^1 + Q_2^3 .
$$
Note that the exponential in factor of $Q_2^1$ is $\exp(i \alpha x + (2 \lambda + \bar \lambda) t)$.
Using Orr Sommerfeld equation we get that
$$
\psi_c = \psi_c^{-3} + \psi_c^{-1} + \psi_c^1 + \psi_c^3,
$$
where in particular
\beq \label{psiqdef}
Orr_{2 \lambda + \bar \lambda,\alpha,\nu}(\psi_c^1) =  {i \over \alpha } Q_2^1.
\eeq
Again $\psi_c^1(0) = \partial_y \psi_c^1(0) = 0$.
Hopefully we do not need to compute explicitly $\psi_c^1$ by using the adjoint of Orr Sommerfeld equation.


\subsection{Adjoint Orr Sommerfeld operator }


We will consider the classical $L^2$ product between two stream functions $\psi_1$ and $\psi_2$, namely
$$
(\psi_1,\psi_2)_\psi = \int \psi_1 \bar \psi_2 dx .
$$
For this scalar product, the adjoint of Orr Sommerfeld equation is
\beq \label{Orrad}
Orr^t_{c,\alpha,\nu}(\psi) :=  (\partial_y^2 - \alpha^2) (U_s - c)   \psi 
- U_s''  \psi  
- { \nu \over i \alpha}   (\partial_y^2 - \alpha^2)^2 \psi,
\eeq
with boundary conditions $\psi(0) = \partial_y \psi(0) = 0$.

We know that the spectrum of Orr Sommerfeld operator and its adjoint are the same. 
Let $\psi_{lin}^t$ be an eigenvector of $Orr_{\lambda,\alpha,\nu}^t$ with corresponding eigenvalue $c$.
Note that
$$
Orr_{2 \lambda + \bar \lambda, \alpha, \nu} = Orr_{\lambda,\alpha,\nu} + \hat c (\partial_y^2 - \alpha^2) 
$$
where
$$
\hat c = c + {2 \lambda + \bar \lambda \over i \alpha} .
$$
Then (\ref{psiqdef}) can be rewritten as
$$
Orr_{\lambda,\alpha,\nu }  \psi_c^1 + \hat c (\partial_y^2 - \alpha^2) \psi_c^1 =  {i \over \alpha } Q_2^1
$$
Taking the scalar product with $\psi_{lin}^t$ we get
$$
\Bigl( Orr_{\lambda,\alpha,\nu }  \psi_c^1, \psi_{lin}^t \Bigr) 
+ \hat c \int  (\partial_y^2 - \alpha^2) \psi_c^1 \bar \psi_{lin}^t =  {i \over \alpha } \int Q_2^1 \bar \psi_{lin}^t .
$$
As $Orr^t_{\lambda,\alpha,c} \psi_{lin}^t = 0$, this gives
\beq \label{useadjoint}
\int \nabla_{\alpha} \psi_c^1 \cdot \nabla_{\alpha} \bar \psi_{lin}^t = - {i \over \alpha  \hat c} \int Q_2^1 \bar \psi_{lin}^t 
\eeq
where $\nabla_\alpha = (\partial_y,i \alpha)$.
Note that $\nabla_\alpha \psi$ is the velocity associated to the stream function $\psi$.
We define the scalar product
$$
(\psi_1,\psi_2)_v = \int \nabla_\alpha \psi_1 \cdot \nabla_\alpha \bar \psi_2 = \int \omega_1 \bar \psi_2 .
$$


\subsection{Computation of $A$}


Let us focus on the mode $+ \alpha$. The approximate solution on this mode is
$$
\psi^{app}_\alpha(t) = \nu^N \psi_{lin} e^{\Re \lambda t + i \Im \lambda t} 
+ \nu^{3 N} \psi_c^1 e^{ 3 \Re \lambda t + i \Im \lambda t} + O(\nu^{5N} e^{ 5 \Re \lambda t}).
$$
Let
$$
\phi(t) =  (\psi^{app}_\alpha , \psi_{lin}^t)_v.
$$ 
We expand the solution $\phi$ of (\ref{model}), which gives
$$
\phi(t) = \nu^N \phi_0 e^{\lambda t}
+ {A \nu^{3N} | \phi_0 |^2 \phi_0 \over 2 \lambda} e^{3 \lambda t} + O(\nu^{5 N} e^{5 \lambda t}).
$$ 
It remains now to compute $A$ by identifying the various terms.
First
$$
\phi_0 = (\psi_{lin},\psi_{lin}^t)_v.
$$
Moreover,
\beq \label{A}
A { | \phi_0 |^2 \phi_0 \over 2 \lambda} =  (\psi_c^1,\psi_{lin}^t)_v,
\eeq
where the right hand side is given by (\ref{useadjoint}).


 \section{Numerical solutions of Orr Sommerfeld \label{part3}}


The computations presented in the previous section can not be done explicitly. 
Numerically they are also difficult since Orr Sommerfeld
equations are very singular near $y = y_c$, where an important part of the dynamics takes place. We therefore need to mix a precise
asymptotic analysis with numerical computations in order to be able to study $A$ as $\nu \to 0$.
This is the objective of the current section.
 
 
 \subsection{Numerical setting}
 
 
 We have to describe the various functions away from the boundary, near the boundary and in the critical layer.
 
 \begin{itemize}
 
 \item Away from the boundary, namely for $\sigma \le y \le Y_0$ we use a grid, with step $h$.
 We approximate Orr Sommerfeld equations by Rayleigh equations, thus neglecting the diffusion of the vorticity in
 this area and we numerically solve Rayleigh equations using an Euler scheme.
 
 \item Close to the boundary we first invert Rayleigh equations by looking for solutions under the form
 \beq \label{formseries}
\psi_{Ray} =  \sum_{n \ge 1} d_n Y^n \log Y 
+ \sum_{n \ge 0} e_n Y^n
 \eeq
 where $Y = y - y_c$.
Doing this we make an error $Diff(\psi_{Ray})$, which is very large and which will be corrected in the next step.
Note that the corresponding vertical velocity has the same form, whereas the corresponding horizontal velocity
has an additional term in $\log Y$ and the vorticity has terms in $\log Y$ and in  $Y^{-1}$.

\item In the critical layer, for $0 \le y \le \Theta \nu^{1/4}$ (with some large $\Theta$), 
we use a grid with a step $h_c$ of order $\nu^{1/4}$
and approximate Orr Sommerfeld equations by ${\cal A}$, with source term $Diff(\psi_{Ray})$.
This leads to an error $ErrAiry({\cal A}^{-1} Diff(\psi_{Ray}))$.

\end{itemize}

Far away from the boundary, namely for $y > Y_0$, we note that the vorticity decays exponentially fast,
whereas the stream function has a slow decay, like $e^{-\alpha y}$.
We will also match $\psi_{Ray}$ and its first derivative at $y = \sigma$, and take care of the boundary conditions at $y = 0$.

Each function of the construction is thus described by its numerical values on two meshes, one ranging from $\sigma$
to $Y_0$ with step $h$, one ranging from $0$ to $\Theta \nu^{1/4}$ with step $h_c$, by three series, one in $y^n \log y$, 
one in $y^n \log^2 y$, and one in $y^n$, used for $0 \le y \le \sigma$.

 
 \subsection{Resolution of the Rayleigh equation \label{series}}
 
 
 
 \subsubsection{Construction of $\psi_{\pm,0}$}
 

We first study the inverse of Rayleigh's operator for $\alpha = 0$. We refer to \cite{GN} for a detailed analysis.
In this case $Ray$ reduces to
  $$
 Ray_0(\psi) = (U_s - c) \partial_y^2 \psi - U_s'' \psi.
 $$
 In particular
 \beq \label{psimoins0}
  \psi_{-,0}(y) = U_s(y) - c
 \eeq
 is an explicit solution of this limiting operator.
 An independent solution $\psi_{+,0}$ is explicitely given by 
 \beq \label{psiplus0}
 \psi_{+,0}(y) =  (U_s(y) - c) \int_0^y{1 \over (U_s(z) - c)^2} dz .
 \eeq
 This other solution behaves linearly at infinity and has a $(y - y_c) \log(y - y_c)$ singularity at $y = y_c$.
 However it is difficult to handle explicit computations using this integral form. Near $y_c$ it is better to look for $\psi_{+,0}$ 
under the form
\beq \label{psiplusnear}
\psi_{+,0}(y) = P(y - y_c) + \log(y - y_c) Q(y -y_c)
\eeq
where $P$ and $Q$ are holomorphic functions near $0$, of the form
$P(Y) = \sum_{n \ge 0} a_n Y^n$ and $Q(Y) = \sum_{n \ge 1} b_n Y^n$.
 Let
$$
Y = y - y_c
$$
and let
$$
  U_s(y) - c = \sum_{n \ge 1} c_n {Y^n \over n! }.
$$
Inserting $U_s$ and $\psi_{+,0}$ in $Ray_0$ we get a series in $Y^n$ and a series in $Y^n \log Y$, which allows to 
compute the various coefficients $a_n$ and $b_n$ by induction, and hence to compute $\psi_{+,0}(y)$ for
$0 \le y \le \sigma$, where $\sigma > 0$ is small enough.
More precisely we have 
$$
\partial_Y^2 \psi_{+,0} = \sum_{n \ge 2} n (n-1) b_n Y^{n-2} \log Y
- \sum_{n \ge 1} b_n Y^{n-2} 
$$
$$
\quad+ 2 \sum_{n \ge 1} n b_n Y^{n-2} + \sum_{n \ge 2} n (n-1) a_n Y^{n-2} .
$$
Writing $(U_s - c) \partial_Y^2 \psi_{+,0} = U_s'' \psi_{+,0}$, we first identify the terms in $Y^0$, leading to
\beq \label{inductionb1}
U_s'(y_c) b_1 = U_s''(y_c) a_0 .
\eeq
Choosing $a_0 = 1$ we thus have 
$b_1 = U_s''(y_c) / U_s'(y_c) $. 
For $n \ge 2$ we then identify the terms in $Y^{n-2} \log Y$ to get $b_n$, defined by induction for $n \ge 3$ by
\beq \label{inductionbn-change}
\begin{split}
&(n-1) (n-2) c_1 b_{n-1} + \sum_{p=2}^{n-2}\frac{ (n-p) (n-p-1) c_p b_{n-p}}{p !}\\
& = \sum_{0 \le p \le n - 3} \frac{(p+2) (p+1) b_{n-2-p} c_{p+2}}{(p+2)!} .
\end{split}
\eeq
Now the coefficients $a_n$ for $n \ge 2$ are defined by induction by identifying terms in $Y^{n-2}$, which leads to
\beq \label{inductionan-change}
\begin{split}&\sum_{1 \le p \le n-1} \frac{c_p (2n - 2p - 1) b_{n-p}}{p!}
	+ \sum_{0 \le p \le n-3} \frac{c_{n-2-p} (p+2) (p+1) a_{p+2}}{(n-2-p)!}  \\
	&\quad = \sum_{0 \le p \le n-2} \frac{a_{n-2-p} (p+2) (p+1) c_{p+2} }{(p+2)!}.
\end{split}
\eeq
This ends the construction of $\psi_{+,0}$ for $0 \le y \le \sigma$.
For $y \ge \sigma$ to get $\psi_{+,0}$ we numerically integrate Rayleigh's equation, 
using the boundary conditions $\psi_{+,0}(\sigma)$
and $\psi_{+,0}'(\sigma)$ given by the previous series.


\subsubsection{Inversion of $Ray_0$}


Let us now turn to the resolution of 
$$
Ray_0(\psi) = \phi
$$
where $\phi$ is a given source term, given on a grid for $y \ge \sigma$ and by a series for $0 \le y \le \sigma$.
For $y \ge \sigma$ we use numerical integration. 
For $y \le \sigma$, if $\phi$ is of the form $P(Y) + Q(Y) \log Y$, then $\psi$ is also of the same form, and the
coefficients of the solution may be explicitly computed by induction. This gives an accurate evaluation
of the solution, even near the singularity $y = y_c$, up to a linear combination of $\psi_{\pm,0}$.
More precisely, for $0 \le y \le \sigma$ let us write $\phi$ under  the form
$$
\phi = \sum_{n \ge 1} d_n Y^n \log Y 
+ \sum_{n \ge 0} e_n Y^n.
$$
We look for $\psi$ under a similar form. We get
\beq \label{inductionb1-source-term}
U_s'(y_c) b_1 = U_s''(y_c) a_0+e_0,
\eeq
\beq \label{inductionbn-source-term}
\begin{split}
&(n-1) (n-2) c_1 b_{n-1} + \sum_{p=2}^{n-2} \frac{(n-p) (n-p-1) c_p b_{n-p} }{p!}\\&
= \sum_{0 \le p \le n - 3}\frac{ (p+2) (p+1) b_{n-2-p} c_{p+2}}{(p+2)!} +d_{n-2},
\end{split}
\eeq
\beq \label{inductionan-source-term}
\begin{split}&\sum_{1 \le p \le n-1} \frac{c_p (2n - 2p - 1) b_{n-p}}{p!}
	+ \sum_{0 \le p \le n-3} \frac{c_{n-2-p} (p+2) (p+1) a_{p+2}}{(n-2-p)!}  \\
	&\quad = \sum_{0 \le p \le n-2} \frac{a_{n-2-p} (p+2) (p+1) c_{p+2} }{(p+2)!}+e_{n-2}.
\end{split}
\eeq

 
 \subsubsection{Rayleigh equation for small $\alpha$}
 

Let us now focus on the case $\alpha > 0$ and study the corresponding Rayleigh operator
which can be rewritten as
\beq \label{Rayleighalpha}
Ray_\alpha(\psi) = (U_s - c + \eps \alpha^2)  \partial_y^2 \psi - \psi \Bigl[ U_s'' + \alpha^2 (U_s - c + \eps \alpha^2) \Bigr].
\eeq
Note that as $U_s''$ converges exponentially fast to $0$, this introduces another asymptotic regime, namely  
$y \gg 1$, with corresponding characteristic scale $\alpha^{-1}$.
We already have an asymptotic expansion of $\psi_{-,\alpha}$, namely
$$
\psi_{-,\alpha}(y) = (U_s(y) - \tilde c) e^{-\alpha y} + O(\alpha),
$$
therefore we can directly evaluate it.

We now turn to $\psi_{+,\alpha}$. For $0 \le y \le \sigma$, we approximate $\psi_{+,\alpha}$ by $\psi_{+,0}$ 
and use the asymptotic expansion obtained in the previous paragraph.
Starting from $y = \sigma$ we then use numerical integration to compute $\psi_{+,\alpha}$ between $\sigma$ and $Y_0$
where $Y_0$ is arbitrarily large.
Note that for large $Y_0$, $\partial_y \psi_{-,\alpha}(Y_0) = O(\alpha)$ whereas 
$\partial_y \psi_{+,\alpha}(Y_0)$ is of order $O(1)$.

We now turn to the resolution of  
$$
Ray_\alpha(\psi) = \phi.
$$ 
Both $\phi$ and $\psi$ are defined by their series between $0$ and $\sigma$,
and by their values on a grid of step $h$ between $\sigma$ and $Y_0$.
We split the computations in two areas:

\begin{itemize}

\item $0 \le y \le \sigma$, where we look for solutions as entire series of the form
(\ref{formseries}). The coefficients of these series are computed by induction. Note that the location $y_c$ of
the singularity slightly changes with $\alpha$. The induction relations are detailed below.

\item $\sigma \le y \le Y_0$, where we directly solve numerically Rayleigh equation, which is regular.

\end{itemize}

Let us detail the computations for $0 \le y \le \sigma$. We have
\beq \label{inductionb1-source-term-alpha}
U_s'(y_c) b_1 = U_s''(y_c) a_0+e_0,
\eeq
\beq \label{inductionbn-source-term-alpha}
\begin{split}
	&(n-1) (n-2) \tilde{c}_1 b_{n-1} + \sum_{p=2}^{n-2} \frac{(n-p) (n-p-1) \tilde{c}_p b_{n-p} }{p!}\\&
	= \sum_{0 \le p \le n - 3}\frac{ (p+2) (p+1) b_{n-2-p} \tilde{c}_{p+2}}{(p+2)!} +\sum_{1 \le p \le n - 3}\frac{\alpha^2\tilde{c}_pb_{n-2-p}}{p!}+d_{n-2},
\end{split}
\eeq
\beq \label{inductionan-source-term-pha}
\begin{split}&\sum_{1 \le p \le n-1} \frac{\tilde{c}_p (2n - 2p - 1) b_{n-p}}{p!}
	+ \sum_{0 \le p \le n-3} \frac{\tilde{c}_{n-2-p} (p+2) (p+1) a_{p+2}}{(n-2-p)!}  \\
	&\quad = \sum_{0 \le p \le n-2} \frac{a_{n-2-p} (p+2) (p+1) \tilde{c}_{p+2} }{(p+2)!}+\sum_{1 \le p \le n - 2}\frac{\alpha^2\tilde{c}_pa_{n-2-p}}{p!}+e_{n-2}.
\end{split}
\eeq


\subsection{Airy equation}



\subsubsection{Construction of the fast mode $\psi_{f,-}$}


In this section we recall the construction of a fast decaying mode $\psi_{f,-}$ of ${\cal A}$. The details may be found in \cite{GN}.
Expanding $U$ near $y_c$ at first order, the modified Airy operator may be approximated by
\beq \label{order2l}
- \eps \partial_y^2 \psi + U_s'(y_c) (y - y_c) \psi = 0 ,
\eeq
which is the classical Airy equation. Let us assume that $\Re U_s'(y_c) > 0$, the opposite case being similar.
A first solution to \eqref{order2l} is given by  
\begin{equation}\label{def-AAA}
A(y) := Ai ( \gamma (y - y_c) ),
\end{equation}
where $Ai$ is the classical Airy function, solution of $Ai'' = x Ai$, and where
$\eps \gamma^3 =  U_s'(y_c)$, namely
$$
\gamma = \Bigl(  {i \alpha U_s'(y_c) \over \nu} \Bigr)^{1/3} .
$$
Note that since $\alpha$ is of order $\nu^{1/4}$, $\gamma$ is of  order $\nu^{-1/4}$ and that 
$$
\arg(\gamma) =  \frac{\pi }{ 6} + O(\nu^{-1/4}).
$$
Another independent solution to \eqref{order2l} is given by $Ci (\gamma (y - y_c))$ where
$$
Ci =  i \pi (Ai + i Bi),
$$
with $Bi(\cdot)$ being the other classical Airy function.

We recall that 
$$
Ai(z) =\sum_{0 \le n \le \infty} \Big[3^{-2/3}\frac{z^{3n}}{  n!\, 9^n\,\Gamma(n+2/3)}  
- 3^{-4/3}\frac{z^{3n+1}} {n!\, 9^n\,\Gamma(n+4/3)}\Big] .
$$
Similarly,
$$
Ai(1,z) =\sum_{0 \le n \le \infty} \Big[\frac{3^{-2/3}z^{3n+1}}{ (3n+1) n!\, 9^n\,\Gamma(n+2/3)}  
- \frac{3^{-4/3}z^{3n+2}} {(3n+2)n!\, 9^n\,\Gamma(n+4/3)}\Big] - 1/3
$$
and
\begin{equation*}
\begin{split}
Ai(2,z) &= \sum_{0 \le n \le \infty} \Big[\frac{3^{-2/3}z^{3n+2}}{ (3n+1)(3n+2) n!\, 9^n\,\Gamma(n+2/3)}  
\\&\quad- \frac{3^{-4/3}z^{3n+3}} {(3n+2)(3n+3)n!\, 9^n\,\Gamma(n+4/3)}\Big]- \frac{z}{3} + 0.25881938,
t\end{split}
\end{equation*}
where the numerical constants insure that both $Ai(1,z)$ and $Ai(2,z)$ go to $0$ as $z$ goes to $+ \infty$.


\subsubsection{Numerical resolution of Airy's equation}


We have to numerically invert ${\cal A} \phi = Airy \, (\partial_y^2 - \alpha^2) \phi = \psi$. We define $\psi$ on a grid
of step $h_c \sim \nu^{1/4}$ near the critical layer. We first numerically solve $Airy \phi_1 = \psi$ assuming that $\phi_1$ together
with its first derivative are negligible away from the critical layer.
Then we numerically solve $(\partial_y^2 - \alpha^2) \phi = \phi_1$ again assuming that $\phi$ and its first derivative are negligible
away from the critical layer.

To solve for $\phi_1$ we use the Green function approach, which allows a good control at infinity.
Let us now detail the Green function $G(x,y)$ of (\ref{order2l}) which by definition is a solution of
\beq \label{order2Green}
- \eps \partial_y^2 G + U_s'(y_c) (y - y_c) G = \delta_x.
\eeq
We have for $y < x$,
$$
G(x,y) = - {Ci(\gamma (y - y_c)) Ai(\gamma (x - y_c)) \over \eps \gamma W}
$$
and for $y > x$
$$
G(x,y) = - {Ai(\gamma(y - y_c)) Ci(\gamma (x - y_c)) \over \eps \gamma W},
$$
where $W$ is the constant Wronskian
$$
W = Ai' Ci - Ci' Ai = 1.
$$
The solution $\phi_1$ to $Airy(\phi_1) = \psi$ is then simply given by
$$
\phi_1(y) = \int G(x,y) \psi(x) dy .
$$

\subsection{Accurate solver for Orr Sommerfeld equations \label{accurate}}


We now focus on the resolution of 
\beq \label{OSSres}
OS(\psi) = f
\eeq
where $f$ is given, of the form $F_{Ray} + F_b$ first without taking care of boundary conditions, then  taking care of them.
Note that $F_{Ray}$ is defined on the "Rayleigh" grid for $\sigma \le y \le Y_0$ and by two series for $0 \le y \le \sigma$,
and $F_b$ is defined on the "Airy" grid, for $0 \le y \le A \gamma^{-1}$ for some large $A$.
We look for $\psi$ of the form
$\psi = \psi_{Ray} + \psi_b$, also defined by two series and on two grids.

\begin{itemize}

\item First we solve the "series" part and compute $\psi_1$, defined by two series between $0$ and $\sigma$, one
in $Y^n$ and the other in $Y^n \log Y$, such that the series of $Ray(\psi_1)$ and of $f_{Ray}$ match.

\item We numerically compute $\psi_1$ for $\sigma \le y \le Y_0$ by integrating Rayleigh equation with source term $f_{Ray}$, 
starting from the boundary conditions $\psi_1(\sigma)$ and $\psi_1'(\sigma)$.

\item For $y > Y_0$, $\psi_1$ is, as a first approximation, a combination of the two independent solutions
$\psi_{-,\alpha}$ and $\psi_{+,\alpha}$, namely $\psi_1 \sim C_- \psi_{-,\alpha} + C_+ \psi_{+,\alpha}$. 
Note that $\partial_y \psi_{-,\alpha}(Y_0) = O(\alpha)$ whereas $\partial_y \psi_{+,\alpha}(Y_0) = O(1)$, 
hence as a first approximation, $C_+ = \psi_1'(Y_0)/ \psi_{+,\alpha}'(Y_0)$.
We thus substract $C_+ \psi_{+,\alpha}$ to $\psi_1$ to ensure a good behavior at infinity.

\item We compute $G_b = Diff(\psi_1 - C_+ \psi_{+,\alpha})$, using series, for $0 \le y \le \sigma$.

\item We invert ${\cal A} \psi_2 = F_b - G_b$, numerically on the "Airy" grid, using Green functions.

\item We adjust the boundary conditions at $y = 0$, by adding a suitable linear combination of $Ai(2,\gamma(y-y_c))$ 
and $\psi_{-,\alpha}(y)$,
which requires to invert the array
$$
A(\alpha) = \left( \begin{array}{cc} 
\psi_{-,\alpha}(0) & Ai(2, - \gamma  y_c ) \cr
\partial_y \psi_{-,\alpha}(0) & \gamma Ai(1,-  \gamma  y_c) \cr 
\end{array} \right) ,
$$
which is singular at $\alpha_0$. Note that $A(2 \alpha_0)^{-1}$ is of size $O(\nu^{-1/4})$.

\end{itemize}


\section{Linear modes  \label{linear}}


This section is devoted to the construction of unstable modes for the Orr Sommerfeld equation, together
with the the corresponding modes of its adjoint.


\subsection{Dispersion relation}


In \cite{GN} it is proved that there exists four independent solutions to $Orr = 0$, two of them decaying at infinity, one "slowly",
called $\psi_{s,-}$ and one "rapidly" called $\psi_{f,-}$.
Moreover
\beq \label{psi1boundary}
\psi_{s,-}(0)= U_s(0) - c + \alpha {U_+^2 \over U_s'(0)} + O(\alpha c),
\eeq
\beq \label{psi1derivboundary}
\partial_y \psi_{s,-}(0) = U_s'(0) + O(\alpha ),
\eeq
and
\beq \label{psif}
\psi_{f,-}(0) = Ai(2, - \gamma y_c), \qquad
\partial_y \psi_{f,-}(0) =  \gamma Ai(1,-\gamma y_c) .
\eeq
An eigenmode of Orr Sommerfeld equation is a combination of these two particular solutions which goes to $0$ at infinity
and which vanishes, together with its first derivative, at $y  = 0$.
There should exist constants $a$ and $b$ such that  $(a,b) \ne (0,0)$,
$$
a \psi_{f,-} (0) + b \psi_{s,-} (0) = 0
$$
and
$$
a \partial_y  \psi_{f,-} (0) + b \partial_y \psi_{s,-} (0) = 0.
$$
The dispersion relation is therefore
\beq \label{disper}
{   \psi_{f,-} (0) \over \partial_y  \psi_{f,-} (0)}
= {  \psi_{s,-} (0) \over \partial_y  \psi_{s,-} (0)}
\eeq
or, using (\ref{psi1boundary}), (\ref{psi1derivboundary}) and  (\ref{psif}),
\beq \label{disper2}
\alpha {U_+^2 \over U_s'(0)^2}  - {c \over U_s'(0)} = \gamma^{-1} {Ai(2, - \gamma y_c) \over Ai(1,-\gamma  y_c)} 
+O(\alpha^2).
\eeq
We will focus on the particular case where $\alpha$ and $c$ are both of order $\nu^{1/4}$. It turns out 
that this is an area where instabilities occur, and we conjecture that this is the region where the most
unstable instabilities may be found.
Therefore we rescale $\alpha$ and $c$ by $\nu^{1/4}$ and introduce
$$
\alpha = \alpha_0 \nu^{1/4}, \quad 
c = c_0 \nu^{1/4}, \quad
Z = \gamma y_c,
$$
which leads to
\beq \label{disper3}
\alpha_0 {U_+^2 \over U_s'(0)^2}  - {c_0 \over U_s'(0)} = 
 {1 \over (i \alpha_0 U_s'(y_c))^{1/3}}  {Ai(2, - Z) \over Ai(1,- Z)} + O(\nu^{1/4}).
\eeq
Note that as $U_s(y_c) = c$, 
$$
y_c = U_s'(0)^{-1} c + O(c)
$$
and
$$
Z =  \Bigl(i  U_s'(y_c) \Bigr)^{1/3} \alpha_0^{1/3} U_s'(0)^{-1} c_0 + O(\nu^{1/4}) .
$$
Note that the argument of $Z$ equals $\pi / 6$.
We then introduce the following function, called Tietjens function, of the real variable $z$
$$
Ti(z) = {Ai(2,z e^{- 5 i \pi / 6}) \over z e^{- 5 i \pi / 6} Ai(1,z e^{- 5 i \pi / 6})} .
$$
At first order the dispersion relation becomes
\beq \label{dispersionlimit}
\alpha_0 {U_+^2 \over U_s'(0)} = c_0 \Bigl[ 1 - Ti(- Z e^{5 i \pi / 6} ) \Bigr] .
\eeq


\subsection{Description of the linear instability \label{desclinear}}


Let us now detail the linear instability for a given $\alpha$. Its stream function $\psi_{lin}$ is of the form
$$
\psi_{lin} = \psi_{s,-} + a \psi_{f,-} + c.c.,
$$
where we have choosen $b = 1$. We see that $a = O(\nu^{1/4})$, hence
\beq \label{psilin}
\psi_{lin} = U_s(y) - c + \alpha {U_+^2 \over U_s'(0)} + a Ai \Bigl(2, \gamma (y - y_c) \Bigr) +  O(\nu^{1/2} ).
\eeq
Note that $a = O(\nu^{1/4})$, hence at leading order $\psi_{lin} = U_s(y) + O(\nu^{1/4})$.

The corresponding horizontal and vertical velocities $u_{lin}$ and $v_{lin}$ are given by
\beq \label{ulin} 
v_{lin} = - i \alpha \psi_{lin} = O(\nu^{1/4}), \qquad u_{lin} = \partial_y \psi_{lin}
\eeq
which at leading order equal
$$
v_{lin} = - i \alpha U_s(y) + O(\nu^{1/2}),
$$
and
\beq \label{vlin}
u_{lin} = \partial_y \psi_{lin} = U_s'(y) + \gamma a Ai \Bigl(1, \gamma (y - y_c) \Bigr) + O(\nu^{1/4}) + c.c..
\eeq
Note that $\gamma a = O(1)$, hence the second term in the right hand side of (\ref{vlin}) is of order $O(1)$ in 
the critical layer. The first term in (\ref{vlin}) may be seen as a "displacement velocity", corresponding
to a translation of $U_s$. The second term is of order $O(1)$ and located in the boundary layer, namely
within a distance $O(\nu^{1/4})$ to the boundary.

The corresponding vorticity is
\beq \label{omegalin}
\omega_{lin} = \omega_{int} + \omega_{bl} = -U_s''(y) - \gamma^2 a Ai \Bigl( \gamma(y - y_c) \Bigr) + \cdots .
\eeq
In particular the vorticity is large, of order $O(\nu^{-1/4})$, near the critical layer and is bounded elsewhere.


\subsection{Construction of adjoint modes}


The construction of $\phi^{t,app}_{s,-}$ has been detailed in \cite{BG1}, where the authors have proven that
\beq \label{approxmode}
\phi^{t,app}_{s,-}(y) = e^{-\alpha y} - f_1(y_c) \psi_3(y) - g_1(y) + \cdots
\eeq
where
	$\psi_3(y)$, $f_1(y_c)$ and $g_1(y)$ satisfy
	$$
	Airy (\psi_3)  = (U_s - c + \eps \alpha^2) \psi_3 - \eps \partial_y^2 \psi_3=e^{-\alpha y} ,
	$$
	$$
	f_1(y) = Ray_\alpha^{-1} \Bigl(- 2 \alpha(U_s - c)  U_s' e^{-\alpha y} \Bigr)
	$$
	and
	$$
	g_1(y)= { f_1(y) - f_1(y_c) \over U_s (y) - c}  .
	$$


\section{Detailed computations for exponential profiles \label{part5}}


We fulfill the previous computations in the particular case of an exponential boundary layer profile
\beq \label{expprofile}
U_s(y) = 1 - e^{-\delta y} .
\eeq
All the routines used may be found in the additional material of this article.
We  discuss the construction in the particular case $\nu = 10^{-30}$, $\delta = 1$ and $\alpha_0 = 1.5$.
In this case the unstable eigenvalue, solution of (\ref{dispersionlimit}), is
$\lambda \sim 4.8 \times 10^{-16} + 3.7 \times 10^{-15} i$, which is of order $\nu^{1/2}$, 
and $y_c \sim c$. 

\medskip

Let us describe the corresponding unstable model (see section \ref{desclinear}). 
First $\psi_{lin}$ is explicit in $Ai(2,.)$ and $U_s$ and is given by (\ref{psilin}). 
Note that there is no $\log(Y)$ part. The corresponding
horizontal velocity $u_{lin}$ goes from $0$ at the boundary to $1$ at a distance of order $0.2 \times 10^{-7}$
and then decays exponentially fast, like $U_s$. The vertical velocity $v_{lin}$ goes smoothly from $0$ to
$-1.5 \times 10^{-7}$ at a distance of order $O(1)$. It is in particular of order $10^{-9}$ in the critical layer. 
It then decays very slowly, like $e^{- \alpha y}$.
The vorticity $\omega_{lin}$ is of order $10^{8}$, namely of order $\nu^{-1/4}$ close to the boundary 
(for $0 \le y \lesssim 0.2 \times 10^{-7}$).

\medskip

Let us turn to the quadratic interaction $Q_1$.
It is small, of order $O(\nu^{1/4})$ outside the critical layer where it is of size $O(1)$.
Note that the source of Orr Sommerfeld equation for $\psi_q^+$ is $i Q_1^+ / 2 \alpha$.
This source term is of order $O(\nu^{-1/4})$ in the critical layer and $O(1)$ inside the flow.
Note that there is no $\log(Y)$ term in $Q_1$.

\medskip

We then follow the strategy sketched in section \ref{accurate} to compute $\psi_q^+$.
First $\psi_1$ and $C_+$ are of order $O(1)$, and $\psi_1 - C_+ \psi_{+,\alpha}$ equals
$\sim 0.6$ at $y = 0$. 
Then $G_b$ is localised in the critical layer and is of order $O(1)$. This is coherent, since 
the most singular term in $G_b$ is $\eps b_1 Y^{-3}$ and thus bounded. Then $\psi_2$ 
is of order $O(1)$ in the boundary layer.
Next we take into account the boundary conditions to construct $\psi_q^+$. As $A(2 \alpha_0)^{-1}$
is of order $O(\nu^{-1/4})$, $\psi_q^+$ is of order $O(\nu^{-1/4})$. 
Note also that $\psi_q^+$ contains terms in $Y^n \log(Y)$.
Moreover, $u_q^+$ has a $\log(Y)$ term in addition to terms in $Y^n$ and $Y^n \log(Y)$ with $n \ge 1$, 
and $\omega_q^+$ has a $\log(Y)$ term and a $Y^{-1}$ term.

\medskip

We now turn to the computation of $Q_2^1$. Note that this term has singular terms in $Y^{-2}$,
$Y^{-1}$, $Y^{-1} \log Y$, $\log Y$ and $\log^2 Y$.

\medskip

The adjoint mode is computed using its explicit formula, and
the scalar product with the adjoint mode is evaluated numerically.

\medskip

For $\alpha_0 = 1.5$  we find  $\Re A \sim - 0.27 < 0$, for $\alpha_0 = 2$, $\Re A \sim -0.03$.
Computations for various $\alpha_0$ show that $\Re A$ is always negative.

%
%
%
%
%
%

\subsubsection*{Acknowledgments}  D. Bian is supported by NSFC under the contract 12271032.




\begin{thebibliography}{99}

\bibitem{BG1} D. Bian, E. Grenier: Long waves instabilities, { \it submitted}, $2022$.


\bibitem{Blasius} H. Blasius: Grenzschichten in Fl\"ussigkeiten mit kleiner Reibung, 
{\it Zeitschrift f\"ur Mathematik und Physik}, Band $56$, Heft $1$, $1908$.

\bibitem{Reid} { P. G.  Drazin, W. H. Reid: } {\em Hydrodynamic stability.} Cambridge Monographs on Mechanics and
Applied Mathematics. Cambridge University, Cambridge--New York, 1981.

\bibitem{GVM}
D.~Gerard-Varet, Y.~Maekawa, N.~Masmoudi:
\newblock Gevrey stability of Prandtl expansions for 2D Navier-Stokes flows, 
{\it Duke Math. J.}, 167(13): 2531-2631, 2018.

\bibitem{GN} {E. Grenier, T. Nguyen:} 
\newblock{Green function for linearized Navier-Stokes around a boundary shear layer profile for long wavelengths,} {\it to appear in Annales Henri Poincar\'e.}
 
 \bibitem{GN2} E. Grenier, T. Nguyen: $L^\infty$ intability of Prandtl layers, {\it Annals of PDE}, $2019$. 
 
\bibitem{GGN3}
{  E. Grenier, Y. Guo, and T. Nguyen:}
\newblock Spectral instability of characteristic boundary layer flows,
\newblock {\em Duke Math. J.}, 165(16): 3085--3146, 2016.

\bibitem{Guo} Y. Guo, S. Iyer: Validity of Steady Prandtl Layer Expansions, {\it  to appear in Comm. Pure Appl. Math.}

\bibitem{Haragus} M. Haragus, G. Iooss: Local bifurcations, center manifolds, and normal forms in infinite dimensional
dynamical systems, Springer, 2011.

\bibitem{Mae}
Y.~Maekawa:
\newblock On the inviscid limit problem of the vorticity equations for viscous
  incompressible flows in the half-plane,
\newblock {\em Comm. Pure Appl. Math.}, 67(7): 1045--1128, 2014.



\bibitem{SammartinoCaflisch1}
M.~Sammartino and R.~E. Caflisch:
\newblock Zero viscosity limit for analytic solutions, of the {N}avier-{S}tokes
  equation on a half-space. {I}. {E}xistence for {E}uler and {P}randtl
  equations,
\newblock {\em Comm. Math. Phys.}, 192(2): 433--461, 1998.

\bibitem{SammartinoCaflisch2}
M.~Sammartino and R.~E. Caflisch:
\newblock Zero viscosity limit for analytic solutions of the {N}avier-{S}tokes
  equation on a half-space. {II}. {C}onstruction of the {N}avier-{S}tokes
  solution,
\newblock {\em Comm. Math. Phys.}, 192(2): 463--491, 1998.



\bibitem{Schensted} I.V. Schensted: Contributions to the theory of hydrodynamic stability, PhD thesis, University of Michigan, $1960$.

\bibitem{Schlichting} { H. Schlichting:} {\em Boundary layer theory,} 
Translated by J. Kestin. 4th ed. McGraw--Hill Series in Mechanical Engineering, McGraw--Hill Book Co., Inc., New York, 1960.



\end{thebibliography}
\end{document}